\documentclass[reqno]{amsart}
\usepackage[dvips]{graphicx}

\theoremstyle{plain}
\newtheorem{theorem}{Theorem}

\newtheorem{lemma}{Lemma}
\newtheorem{corollary}{Corollary}
\theoremstyle{definition}

\begin{document}

\title{Finite-gap minimal Lagrangian surfaces in ${\mathbb C}P^2$}
\author{A.E. Mironov}

\address{Sobolev Insitute of Mathematics and Novosibirsk
State University, 630090 Novosibirsk, Russia}
\email{mironov@math.nsc.ru }

\begin{abstract}
In this paper we suggest a method for constructing minimal
Lagrangian immersions of ${\mathbb R}^2$ in ${\mathbb C}P^2$ with
induced diagonal metric in terms of Baker--Akhiezer functions of
algebraic curves.
\end{abstract}

\thanks {This research was partially supported by the Russian
Foundation for Basic Research (grant no. 09-01-92442-KEa) and
grant MK-5430.2008.1 of the President of Russian Federation.}

\maketitle
\section*{Introduction}
We propose new approach for constructing minimal Lagrangian ({\bf
ML}) surfaces in ${\mathbb C}P^2$ in terms of Baker--Akhiezer
functions of algebraic curves. This approach is based on the work
\cite{M1}.

It is well known that if we choose the conformal coordinates on a
{\bf ML}-torus in ${\mathbb C}P^2$ with the induced metric
$ds^2=e^{v(x,y)}(dx^2+dy^2)$, then the function $v(x,y)$ satisfies
the Tzitzeica equation (see \cite{HM}). This equation allows the Lax
representation with a spectral parameter, which was found by
Mikhailov \cite{Mih}. Sharipov \cite{Sh}, using the Lax
representation, constructed the finite-gap solutions of the
Tzitzeica equation, and the solutions expressed in terms of the
theta-functions of the trigonal spectral curves, which allow the
holomorphic involution. The existence of periodic solutions among
quasiperiodic solutions is shown in \cite{Mc}.

In fact, the conformal coordinates are not always suitable for the
description of {\bf ML}-tori in ${\mathbb C}P^2$. To confirm this,
let us consider the following example. Let $K$ denote a cone in
${\mathbb R}^3$, defined by the equation
$$
 mu_1^2+nu_2^2=(m+n)u_3^2,
$$
where $m,n\in{\mathbb Z}, m,n>0.$ Let ${\tilde K}$ denote the
intersection of $K$ with the unit sphere
$$
 u_1^2+u_2^2+u_3^2=1.
$$
We construct the mapping $\psi$ from ${\tilde K}\times S^1$ in
${\mathbb C}P^2$ as a composition $\psi= \mathcal{H}\circ
\tilde{\varphi},$ where
$$
 \tilde{\varphi}:{\tilde K}\times S^1\rightarrow S^5,
 \ {\tilde \varphi}(P)=(u_1e^{\pi imy},u_2e^{\pi iny},u_3e^{\pi i(m+n)y}),
$$
$\mathcal{H}$ --- Hopf projection $\mathcal{H}:S^5\rightarrow
{\mathbb C}P^2$, $P\in{\tilde K}\times S^1$, $y$ is a coordinate on
$S^1$.

The image of $\psi$ is {\bf ML}-torus if the involution
$$
 (v_1,v_2)\rightarrow(v_1\cos(n\pi),v_2\cos(m\pi))
$$
preserve the orientation of the ellipse
$$
 mv_1^2+nv_2^2=m+n
$$
and {\bf ML}-Klein bottle if it doesn't (see \cite{M2},\cite{J}).
This surface can be defined as an image of the composition
$\mathcal{H}\circ\varphi,$ where $\varphi:{\mathbb R}^2\rightarrow
S^5,$
$$
  \varphi(x,y)=
 \left(\frac{\sin(x)\sqrt{m+n}}{\sqrt{2m+n}}e^{\pi imy},\right.
$$
$$
 \left.\frac{\cos(x)\sqrt{m+n}}{\sqrt{m+2n}}e^{\pi i ny},
  \sqrt{\frac{n\cos^2(x)}{m+2n}+\frac{m\sin^2(x)}{2m+n}}e^{-\pi i(m+n)y}\right)
$$
(probably these tori coincide with the tori from \cite{CU}, where these tori are described in
conformal coordinates).
 In coordinates $x,y$ the induced metric has a diagonal form
$$
 ds^2=2e^{v_1}dx^2+2e^{v_2}dy^2,
$$
and one can ascertain that $v_1 \ne v_2$. Thus, on the one hand {\bf
ML}-tori correspond to the periodic solutions of the Tzitzeica
equation (all of these solutions can be expressed in terms of the
theta-function of the spectral curves), and on the other hand this
example demonstrates that there exist coordinates $x,y$, in which
the metric is diagonal and which are more suitable for the
description of this {\bf ML}-torus, as in these coordinates the
tori are described in elementary functions.

In this paper we construct {\bf ML}-mapping of the plane in
${\mathbb C}P^2$ (with induced diagonal metric) by the spectral
data, which are easier than spectral data for the solution of the
Tzitzeica equation. We construct such mapping by the real
algebraic curve, which allows a holomorphic involution.
Specifically, we do not require the spectral curve to be trigonal
(as for the solutions of the Tzitzeica equation).

The main difference of our method from the method of \cite{Sh}
consists in the following. We do not use the Lax representation
with a spectral parameter of the Tzitzeica equation. Instead of
this we construct the explicit mapping
$$
 \varphi:{\mathbb R}^2\rightarrow S^5\subset {\mathbb C}^3,
$$
which satisfies the equations
$$
 <\varphi,\varphi_x>=<\varphi,\varphi_y>=<\varphi_x,\varphi_y>=0,\eqno{(1)}
$$
where $<.,.>$ --- Hermitian product in ${\mathbb C}^3$. A
composition of the mappings $\varphi\circ \mathcal{H}$ gives a
Lagrangian mapping of the plane in ${\mathbb C}P^2$. By means of
the corollary 1 we obtain the minimal mappings.

Note that in this paper we do not discuss the problem of existence
of the smooth periodic solutions. As the spectral curve is
hyperelliptic, the methods of the paper \cite{E} can be used to
prove the existence of the periodic solutions.

In section 1 we get the equations of  {\bf ML}-mapping of a plane
in ${\mathbb C}P^2$ with diagonal metric. In section 2 we remind
the definition of the Baker--Akhiezer function. In section 3 we
prove the main theorem (Theorem 2) and give the example of  {\bf
ML}-sphere, corresponding to the reducible rational spectral
curve.

\section{Equations of Lagrangian surfaces with the diagonal metric }
Let we introduce the following notations
$$
 |\varphi_x|^2=2e^{v_1(x,y)},\ \ |\varphi_y|^2=2e^{v_2(x,y)}.
$$
Then from (1) it follows, that the matrix
$$
 \tilde{\Phi}=
 \left(\varphi,\frac{1}{\sqrt{2}}e^{-\frac{v_1}{2}}\varphi_x,
 \frac{1}{\sqrt{2}}e^{-\frac{v_2}{2}}\varphi_y\right)^{\top}
$$
belongs to the group ${\bf U(3)}$. {\it A Lagrangian angle
$\beta(x,y)$} is a function defined from the equality
$$
 e^{i\beta (x,y)}={\rm det} \tilde{\Phi}.
$$
The Lagrangian angle defines the mean curvature vector $H$ of the
Lagrangian surface
$$
 H=J\nabla\beta,
$$
where $J$ is the complex structure on ${\mathbb C}P^2$.
Consequently, if $\beta=const$, then the surface is minimal.

From the definition of the Lagrangian angle we get
$$
 \Phi=
  \left(
  \begin{array}{ccc}
    \varphi^1 & \varphi^2 &  \varphi^3\\
   \frac{1}{\sqrt{2}}e^{-\frac{v_1}{2}-i\frac{\beta}{2}}\varphi^1_x
   &  \frac{1}{\sqrt{2}}e^{-\frac{v_1}{2}-i\frac{\beta}{2}}\varphi^2_x &
    \frac{1}{\sqrt{2}}e^{-\frac{v_1}{2}-i\frac{\beta}{2}}\varphi^3_x \\
    \frac{1}{\sqrt{2}}e^{-\frac{v_2}{2}-i\frac{\beta}{2}}\varphi^1_y
    &  \frac{1}{\sqrt{2}}e^{-\frac{v_2}{2}-i\frac{\beta}{2}}\varphi^2_y
    &  \frac{1}{\sqrt{2}}e^{-\frac{v_2}{2}-i\frac{\beta}{2}}\varphi^3_y \\
  \end{array}\right)\in{\bf SU(3)}.
$$
The matrix $\Phi$ satisfies the equations
$$
 \Phi_x=A\Phi,\ \Phi_y=B\Phi,\eqno{(2)}
$$
where matrices $A, B\in{\bf su(3)}$ have the form
$$
 A=
  \left(
  \begin{array}{ccc}
   0 & \sqrt{2}e^{\frac{v_1}{2}+i\frac{\beta}{2}} & 0\\
  -\sqrt{2}e^{\frac{v_1}{2}-i\frac{\beta}{2}} & if &
  \frac{1}{2}e^{\frac{v_1}{2}-\frac{v_2}{2}}(2ih-v_{1_{y}}+i\beta_y) \\
  0 & \frac{1}{2}e^{\frac{v_1}{2}-\frac{v_2}{2}}(2ih+v_{1_{y}}+i\beta_y) & -if\\
  \end{array}\right),
$$
$$
 B=
  \left(
  \begin{array}{ccc}
   0 & 0 & \sqrt{2}e^{\frac{v_2}{2}+i\frac{\beta}{2}} \\
  0 & ih & \frac{1}{2}e^{\frac{v_2}{2}-\frac{v_1}{2}}(i\beta_x-2if +v_{2_x})  \\
  -\sqrt{2}e^{\frac{v_2}{2}-i\frac{\beta}{2}} &
   \frac{1}{2}e^{\frac{v_2}{2}-\frac{v_1}{2}}(i\beta_x-2if-v_{2_x})  & -ih\\
  \end{array}\right),
$$
$f(x,y)$ and $h(x,y)$  being some real functions. The equations
(2) implies the following equations
$$
 \varphi_{xx}=\Gamma_{11}^1\varphi_x+\Gamma_{11}^2\varphi_y+b_{11}\varphi,
$$
$$
 \varphi_{xy}=\Gamma_{12}^1\varphi_x+\Gamma_{12}^2\varphi_y+b_{12}\varphi,
$$
$$
 \varphi_{yy}=\Gamma_{22}^1\varphi_x+\Gamma_{22}^2\varphi_y+b_{22}\varphi,
$$
where
$$
 \Gamma_{11}^1=\frac{1}{2}(2if+v_{1_x}+i\beta_x),\ \
 \Gamma_{11}^2=\frac{1}{2}e^{v_1-v_2}(2ih-v_{1_{y}}+i\beta_y),
$$
$$
 \Gamma_{12}^1=\frac{1}{2}(2ih+v_{1_y}+i\beta_y),\ \
 \Gamma_{12}^2=\frac{1}{2}(-2if+v_{2_x}+i\beta_x),
$$
$$
 \Gamma_{22}^1=\frac{1}{2}e^{v_2-v_1}(-2if-v_{2_{x}}+i\beta_x),\
 \Gamma_{22}^2=\frac{1}{2}(-2ih+v_{2_y}+i\beta_y),
$$
$$
 b_{11}=-2e^{v_1},\ b_{12}=0,\ b_{22}=-2e^{v_2}.
$$
From this it follows the key lemma of our construction.

\begin{lemma}
The following equalities hold:
$$
 \Gamma_{11}^1+\Gamma_{12}^2=\frac{1}{2}(v_{1_x}+v_{2_x})+i\beta_x,
$$
$$
 \Gamma_{12}^1+\Gamma_{22}^2=\frac{1}{2}(v_{1_y}+v_{2_y})+i\beta_y.
$$
\end{lemma}

From lemma 1 we  get

\begin{corollary}
If
$$
 {\rm Im}(\Gamma_{11}^1+\Gamma_{12}^2)= {\rm Im}(\Gamma_{12}^1+\Gamma_{22}^2)=0,
$$
then the surface is minimal.
\end{corollary}

The corollary 1 gives us the condition for surface to be minimal
in a diagonal metric.

\section{Baker--Akhiezer function}

In this paragraph we remind the definition of two-point
Baker--Akhiezer function. By means of this function we construct
{\bf ML}-mapping of the plane in ${\mathbb C}P^2$.

Let $\Gamma$ be a Riemann surface of genus $g$ (actually the
following construction can be generalized on singular algebraic
curves over ${\mathbb C}$). Suppose that the divisor
$$
 \gamma=\gamma_1+\dots+\gamma_{g},
$$
is given on $\Gamma$,  $r$, $P_1$, $P_2\in \Gamma$ are fixed
points, and  $k_1^{-1}, k_2^{-1}$ are  local parameters in the
neighborhoods of the points $P_1$ and $P_2$. The {\it two-point
Baker--Akhiezer function}, corresponding to the  {\it spectral
data}
$$
 \{\Gamma,P_1,P_2,k_1,k_2,\gamma,r\},
$$
is a function $\psi(x,y,P), P\in\Gamma$, with the following characteristics:

{\bf 1)} in the neighborhoods of $P_1$ and $P_2$ the function
$\psi$ has essential singularities of the following form:
$$
 \psi=e^{ik_1x}\left(f_1(x,y)+\frac{g_1(x,y)}{ik_1}+\frac{h_1(x,y)}{k_1^2}+\dots\right),
$$
$$
 \psi=e^{ik_2y}\left(f_2(x,y)+\frac{g_2(x,y)}{ik_2}+\frac{h_2(x,y)}{k_2^2}+\dots\right).
$$

{\bf 2)} on $\Gamma\backslash \{P_1,P_2\}$ the function $\psi$
is meromorphic with simple poles on $\gamma$.

{\bf 3)} $\psi(x,y,r)=d$, $d\in{\mathbb C}.$

For the spectral data in general position there is unique
Baker--Akhiezer function.

Let us express the Baker--Akhiezer function explicitly in terms of
theta function of the surface $\Gamma$.

On the surface $\Gamma$, choose a basis of cycles
$$
 a_1,\dots,a_g, \ b_1,\dots,b_g
$$
with the following intersections indices:
$$
 a_i\circ a_j=b_i\circ b_j=0,\ a_i\circ b_j=\delta_{ij}.
$$
By $\omega_1,\dots,\omega_g$ we denote a basis of holomorphic differentials on $\Gamma$
that are normalized by the conditions
$$
 \int_{a_j}\omega_i=\delta_{ij}.
$$
Denote the matrix of $b$-periods of the differentials $\omega_j$ with the components
$$
 B_{ij}=\int_{b_i}\omega_j
$$
by $B$. This matrix is symetric and has a positive definite imaginary part.

The Riemann theta function is defined by the absolutely converging
series
$$
 \theta(z)=\sum_{m\in{\mathbb Z}}e^{\pi i (Bm,m)+2\pi i(m,z)},\
 z=(z_1,\dots,z_g)\in{\mathbb C}^g.
$$
The theta function has the following characeristics:
$$
 \theta(z+m)=\theta(z),\ m\in{\mathbb Z},
$$
$$
 \theta(z+Bm)={\rm exp}(-\pi i (Bm,m)-2\pi i (m,z))\theta(z),\ m\in{\mathbb
 Z}.
$$
Let $X$ denote the Jacobi variety of the surface $\Gamma$:
$$
 X={\mathbb C}^g/ \{{\mathbb Z}^g+B{\mathbb Z}^g\}.
$$
Let $A:\Gamma\rightarrow X$ be an Abel map defined by the formula
$$
 A(P)=\left(\int_{q_0}^P\omega_1,\dots,\int_{q_0}^P\omega_g\right),\
 P\in\Gamma,
$$
$q_0\in\Gamma$ being a fixed point.

For points $\gamma_1,\dots,\gamma_g$ in general position, according to the Riemann theorem, the function
$$
 \theta(z+A(P)),
$$
where $z=K-A(\gamma_1)-\dots-A(\gamma_g)$, has exacly $g$ zeros
$\gamma_1,\dots,\gamma_g$ on $\Gamma$, $K$ is the vector of
Riemann constants.

Let $\Omega^1$ and $\Omega^2$ denote meromorphic differentials on $\Gamma$
with simple poles only at the points $P_1$ and $P_2$, respectively, and normalized by the conditions
$$
 \int_{a_j}\Omega^1=\int_{a_j}\Omega^2=0,\ j=1,\dots,g.
$$
Let $U$ and $V$ denote their vectors of $b$-periods:
$$
 U=\left(\int_{b_1}\Omega^1,\dots,\int_{b_g}\Omega^1\right),\
 V=\left(\int_{b_1}\Omega^2,\dots,\int_{b_g}\Omega^2\right).
$$
Let $\widetilde{\psi}$ denote the function
$$
 \widetilde{\psi}(x,y,P)=\frac{\theta(A(P)+xU+yV+z)}{\theta(A(P)+z)}
 {\rm exp}(2\pi ix\int_{q_0}^P\Omega^1+2\pi iy\int_{P_1}^P\Omega^2).
$$
Then the Baker--Akhiezer function has the form:
$$
 \psi(x,y,P)=\frac{\widetilde{\psi}(x,y,P)}{\widetilde{\psi}(x,y,r)}d.
$$

\section{Main theorem}

Let $\varphi^1,\varphi^2,\varphi^3$ denote the following functions:
$$
 \varphi^i=\alpha_i\psi(x,y,Q_i),
$$
where $Q_1,Q_2,Q_3\in\Gamma$ is an additional set of points and $\alpha_i$
are some constants.

In the following paragraph (theorem 1) we find restrictions of the
spectral data for the vector valued function
$\varphi=(\varphi^1,\varphi^2,\varphi^3)$ to define a Lagrangian
immersion of the plane ${\mathbb R}^2$ in ${\mathbb C}P^2$. The
spectral data in theorem 1 is a modification of the spectral data
for $n$-orthogonal curvilinear coordinate system in ${\mathbb
R}^3$ \cite{Kr}.

\subsection{Lagrangian immersions}
Suppose that the surface $\Gamma$ has an antiholomorphic
involution $\mu$
$$
  \mu:\Gamma\rightarrow\Gamma,
$$
for which the points $P_1$, $P_2$ and $r$ are fixed and
$$
 k_i(\mu(P))=\bar{k}_i(P).
$$
The following theorem holds:

\begin{theorem}
Let $Q_i$ be fixed points of the antiholomorphic involution $\mu$. Suppose that
on $\Gamma$ there exists a meromorpic 1-form $\Omega$ with the following set
of divisors of zeros and poles:
$$
(\Omega)_0=\gamma+\mu\gamma+P_1+P_2,
$$
$$
 (\Omega)_{\infty}=Q_1+Q_2+Q_3+r.
$$
Then the functions $\varphi^i$ satisfy the equations
$$
 \varphi^1\bar{\varphi}^1A_1+
 \varphi^2\bar{\varphi}^2A_2+\varphi^3\bar{\varphi}^3A_3+
 |d|^2\rm{Res}_{r}\Omega=0,
$$
$$
 \varphi^1\bar{\varphi}^1_xA_1+
 \varphi^2\bar{\varphi}^2_xA_2+
 \varphi^3\bar{\varphi}^3_xA_3=0,
$$
$$
 \varphi^1\bar{\varphi}^1_yA_1+
 \varphi^2\bar{\varphi}^2_yA_2+
 \varphi^3\bar{\varphi}^3_yA_3=0,
$$
$$
 \varphi^1_x\bar{\varphi}^1_yA_1+
 \varphi^2_x\bar{\varphi}^2_yA_2+
 \varphi^3_x\bar{\varphi}^3_yA_3=0,
$$
$$
 \varphi^1_x\bar{\varphi}^1_xA_1+
 \varphi^2_x\bar{\varphi}^2_xA_2+
 \varphi^3_x\bar{\varphi}^3_xA_3+|f_1|^2c_1=0,\eqno{(4)}
$$
$$
 \varphi^1_y\bar{\varphi}^1_yA_1+
 \varphi^2_y\bar{\varphi}^2_yA_2+
 \varphi^3_y\bar{\varphi}^3_yA_3+|f_2|^2c_2=0,\eqno{(5)}
$$
where $A_k=\frac{\rm{Res}_{Q_k}\Omega}{|\alpha_k|^2},\ k=1,2,3$,
$c_1, c_2$ are the coefficients of the form $\Omega$ in
the neighborhood of the points $P_1$ and $P_2$:
$$
 \Omega=(c_1w_1+aw_1^2+\dots)dw_1,\ w_1=1/k_1,
$$
$$
 \Omega=(c_2w_2+bw_2^2+\dots)dw_2,\ w_2=1/k_2.
$$
\end{theorem}

The theorem 1 implies
\begin{corollary}
If $\rm{Res}_{Q_i}\Omega>0,$ then for
$$
 \alpha_i=\sqrt{\rm{Res}_{Q_i}\Omega},\
 d=\sqrt{\frac{-1}{\rm{Res}_{r}\Omega}},
$$
the following equalities hold:
$$
   <\varphi,\varphi>=1,\ <\varphi,\varphi_x>=<\varphi,\varphi_y>=<\varphi_x,\varphi_y>=0,
$$
i.e. the mapping $\mathcal{H}\circ \varphi:{\mathbb
R}^2\rightarrow {\mathbb C}P^2$ is Lagrangian, with the induced
metric on $\Sigma$ having a diagonal form
$$
 ds^2=|f_1|^2|c_1|dx^2+|f_2|^2|c_2|dy^2.
$$
\end{corollary}
Further we assume that $\rm{Res}_{Q_i}\Omega>0$ and $f_1\ne 0,
f_2\ne 0.$

 {\bf Proof of Theorem 1}. Consider the 1-form
$\Omega_1=\psi(P)\overline{\psi(\mu (P))}\Omega$. By virtue of the
definition of the involution $\mu$, the function
$\overline{\psi(\mu (P))}$ has the following form in the
neighborhoods of the points $P_1$ and $P_2$:
$$
 \overline{\psi(\mu (P))}=e^{-i k_1x}
 \left(\bar{f}_1(x,y)-\frac{\bar{g}_1(x,y)}{ik_1}+\frac{\bar{h}_1(x,y)}{k_1^2}+\dots\right),
$$
$$
 \overline{\psi(\mu (P))}=e^{-i k_2y}\left(\bar{f}_2(x,y)-
 \frac{\bar{g}_2(x,y)}{ik_2}+\frac{\bar{h}_2(x,y)}{k_2^2}+\dots\right).
$$
Consequently, the form $\Omega_1$ has no essential singularities
at the points $P_1$ and $P_2$. The simple poles $\gamma+\mu\gamma$
of the function $\psi(P)\overline{\psi(\mu P)}$ cancel out the
zeros of the form $\Omega$ at these points. Thus, the form
$\Omega_1$ has only simple poles at the points $Q_1,Q_2,Q_3$ and
$r$ with the residues equal to
$$
 \psi(Q_1)\overline{\psi(Q_1)}{\rm Res}_{Q_1}\Omega=
 \varphi^1\bar{\varphi}^1A_1,\
 \varphi^2\bar{\varphi}^2A_2,\ \varphi^3\bar{\varphi}^3A_3,\
 |d|^2\rm{Res}_{r}\Omega.
$$
Consequently, the sum of these residues is equal to zero, and this proves the
first equality of the theorem 1.

The form $\psi(P)\overline{\psi(\mu (P))_x}\Omega$ has no
essential singularities at the points $P_1$ and $P_2$ either. This
form has only simple poles at the points $Q_1,Q_2$ and $Q_3$ with
the residues equal to
$$
 \varphi^1\bar{\varphi}^1_xA_1,\
 \varphi^2\bar{\varphi}^2_xA_2,\
 \varphi^3\bar{\varphi}^3_xA_3.
$$
The second equality is proven. The proof of the next two equalities is analogous.
It is based on the analysis of the forms
$$
 \psi(P)\overline{\psi(\mu (P))_y}\Omega,\
 \psi(P)_x\overline{\psi(\mu (P))_y}\Omega,
$$
which also have only simple poles at the points $Q_1,Q_2$ è $Q_3$.

In order to prove the last two equalities (4) and (5), consider
the forms
$$
 \psi(P)_x\overline{\psi(\mu (P))_x}\Omega,
\
 \psi(P)_y\overline{\psi(\mu (P))_y}\Omega.
$$
These forms have simple poles at the points $Q_1,Q_2$, $Q_3, P_1$ and $Q_1,Q_2$, $Q_3, P_2$
with the residues
$$
 \varphi^1_x\bar{\varphi}^1_xA_1,\
 \varphi^2_x\bar{\varphi}^2_xA_2,\
 \varphi^3_x\bar{\varphi}^3_xA_3,\ |f_1|^2c_1
$$
and
$$
 \varphi^1_y\bar{\varphi}^1_yA_1,\
 \varphi^2_y\bar{\varphi}^2_yA_2,\
 \varphi^3_y\bar{\varphi}^3_yA_3,\ |f_2|^2c_2.
$$
Theorem 1 is proven.

\subsection{Minimal Lagrangian immersions}
In this subsection we find spectral data such that the mapping
$\varphi$, constructed in the previous  subsection is minimal.

Suppose that the curve $\Gamma$ has the holomorphic involution
$$
 \sigma:\Gamma\rightarrow \Gamma.
$$
Let $\tau$ denote the composition $\mu\circ \sigma.$ The following
lemma holds.

\begin{lemma}
Suppose that the reality conditions fulfilled
$$
 \mu(\gamma)=\gamma, \ \mu(r)=r,\ d\in{\mathbb R}.
$$
Then
$$
 \psi(x,y,\tau(P))=\overline{\psi(x,y,P)}.
$$
\end{lemma}

To prove this standard lemma, it is sufficient to note, that the
function $\overline{\psi(x,y,\tau(P))}$ satisfies conditions {\bf
1)}--{\bf 3)} in the definition of the Baker--Akhiezer function as
the function $\psi(x,y,P)$, consequently, the functions
$\overline{\psi(x,y,\tau(P))}$ and $\psi(x,y,P)$ coincide.

Below we assume that the conditions of Lemma 2 are fulfilled. In
particular, this means that
$$
 \mu(\gamma)=\sigma(\gamma), \ \mu(r)=\sigma(r),
$$
and that the functions $f_i,g_i$, which participate in the
decomposition of $\psi$ in the neighborhoods of the points $P_1$
and $P_2$ are real.

Consider three functions
$$
 F_{11}(x,y,P)=\partial_x^2\psi+\Gamma_{11}^1(x,y)\partial_x\psi+
 \Gamma_{11}^2(x,y)\partial_y\psi+b_{11}(x,y)\psi,
$$
$$
 F_{12}(x,y,P)=\partial_x\partial_y\psi+\Gamma_{12}^1(x,y)\partial_x\psi+
 \Gamma_{12}^2(x,y)\partial_y\psi+b_{12}(x,y)\psi,
$$
$$
 F_{22}(x,y,P)=\partial_y^2\psi+\Gamma_{22}^1(x,y)\partial_x\psi+
 \Gamma_{22}^2(x,y)\partial_y\psi+b_{22}(x,y)\psi.
$$
Choose functions $ \Gamma_{ij}^k(x,y)$ and $b_{ij}(x,y)$ such that
$$
 F_{11}(x,y,Q_i)=F_{12}(x,y,Q_i)=F_{22}(x,y,Q_i)=0,\ i=1,2,3.
$$
The following lemma holds

\begin{lemma} The following equalities hold:
$$
 \Gamma_{11}^1(x,y)=-\frac{ia}{c_1}-\frac{f_{1_x}}{f_1},
$$
$$
 \Gamma_{12}^1(x,y)=-\frac{f_{1_y}}{f_1},\ \Gamma_{12}^2(x,y)=-\frac{f_{2_x}}{f_2},
$$
$$
 \Gamma_{22}^2(x,y)=-\frac{ib}{c_2}-\frac{f_{2_y}}{f_2}.
$$
\end{lemma}
{\bf Proof of lemma 3.} Consider the form
$$
 \omega=F_{11}(P)\psi(\sigma(P))_x\Omega.
$$
The form $\omega$ has no essential singularities at the points
$P_1$ and $P_2$. The form $\omega$ has only a pole of the second
order at $P_1$
$$
 {\rm Res}_{P_1}\omega=f_1((ia+c_1\Gamma_{11}^1)f_1+c_1f_{1_x})=0.
$$
This yields the formula for $\Gamma_{11}^1.$ Similarly to find
other coefficients it is necessary to consider the forms
$$
F_{12}(P)\psi(\sigma(P))_x\Omega,\ F_{12}(P)\psi(\sigma(P))_y\Omega,\ F_{22}(P)\psi(\sigma(P))_y\Omega.
$$
Lemma 3 is proven.

This lemma implies
$$
 \Gamma_{11}^1+\Gamma_{12}^2=-\frac{ia}{c_1}-\frac{f_{1_x}}{f_1}-\frac{f_{2_x}}{f_2},
$$
$$
 \Gamma_{12}^1+\Gamma_{22}^2=-\frac{ib}{c_2}-\frac{f_{1_y}}{f_1}-\frac{f_{2_y}}{f_2}.
$$
Suppose that $a=b=0$. Then as $f_1$ and $f_2$ are real functions,
so by Corollary 1 the surface is minimal. Thus the main theorem
holds

\begin{theorem}
Suppose that the spectral curve $\Gamma$ has the antiholomorphic
involu\-tion
$$
 \mu:\Gamma\rightarrow\Gamma
$$
with fixed points $Q_1,Q_2,Q_3,P_1,P_2$ è $r$ and meromorfic 1-form $\Omega$ with the following
divisors of zeros and poles
$$
 (\Omega)_0=\gamma+\mu\gamma+P_1+P_2,\ (\Omega)_{\infty}=Q_1+Q_2+Q_3+r,
$$
and ${\rm Res}_{Q_i}\Omega>0$. Then the mapping $\mathcal{H}\circ
\varphi$, where $\varphi=(\varphi_1,\varphi_2,\varphi_3)$ gives
Lagrangian mapping of the plane in ${\mathbb C}P^2$.

Besides let the spectral curve $\Gamma$ has the holomorphic
involution
$$
 \sigma:\Gamma\rightarrow \Gamma
$$
such that
$$
 \mu(\gamma)=\sigma(\gamma),\ \tau(r)=r
$$
and $d\in{\mathbb R}$ and suppose, that the form $\Omega$ has the following decomposition in the neighborhood
of the points $P_1$ è $P_2$
$$
 \Omega=(c_1w_1+d_1w_1^3+\dots)dw_1,\ w_1=1/k_1,\eqno(6)
$$
$$
 \Omega=(c_2w_2+d_2w_2^3+\dots)dw_2,\ w_2=1/k_2,\eqno(7)
$$
then the mapping is minimal.
\end{theorem}

\subsection{Examples} In this paragraph we demonstrate the example
of theorem 2, when the spectral curve $\Gamma$ is reducible and
consists of irreducible components $\Gamma_i$, which are
isomorphic ${\mathbb C}P^1$. In this case the theorem 2 is also
valid, but in definition of Baker--Akhiezer function the genus
need to be changed on arithmetical genus, and in the formulation
of theorem 2 the differential need to be changed on the
differential, which satisfies the condition of regularity at the
points of intersection of different components.

The  regular differential on $\Gamma$ is defined by meromorphic
1-forms $\Omega_j$ on $\Gamma_j$ with simple poles. The poles of
the forms $\Omega_j$ are allowed just in the points of components'
intersections. And the conditions of regularity must be fulfilled:
if the curves $\Gamma_i$ and $\Gamma_j$ intersect at the point
$P$, then
$$
 {\rm Res}_P\Omega_1+{\rm Res}_P\Omega_2=0.
$$
The arithmetical genus of the curve $\Gamma$ is called the
dimension of the space of the regular differentials. A number of
poles $\gamma_i$ in the definition of Baker--Akhiezer function
must coincide with the arithmetical genus of the curve $\Gamma$
(see \cite{MT}).

Let the curve $\Gamma$ consists of two components $\Gamma_1$ and
$\Gamma_2$, which intersect at two points.

\vskip15mm

\begin{picture}(170,125)(-100,-80)
\qbezier(-10,0)(-10,30)(40,30)
\qbezier(40,30)(90,30)(90,0)
\qbezier(90,0)(90,-30)(40,-30)
\qbezier(40,-30)(-10,-30)(-10,0)
\put(35,33){\shortstack{$\Gamma_1$}}

\qbezier(60,0)(60,30)(110,30)
\qbezier(110,30)(160,30)(160,0)
\qbezier(160,0)(160,-30)(110,-30)
\qbezier(110,-30)(60,-30)(60,0)
\put(105,33){\shortstack{$\Gamma_2$}}

\put(-10,0){\circle*{3}}
\put(90,0){\circle*{3}}
\put(60,0){\circle*{3}}

\put(75,24){\circle*{3}}
\put(75,-24){\circle*{3}}

\put(110,-30){\circle*{3}}
\put(134,-28){\circle*{3}}
\put(155,-17){\circle*{3}}

\put(-25,0){\shortstack{$P_1$}}
\put(105,-44){\shortstack{$Q_1$}}
\put(132,-43){\shortstack{$Q_2$}}

\put(159,-29){\shortstack{$Q_3$}} \put(95,0){\shortstack{$r$}}
\put(45,0){\shortstack{$P_2$}}
\put(63,20){\shortstack{\small{$a$}}}
\put(54,-26){\shortstack{\small{$-a$}}}
\put(85,18){\shortstack{\small{$b$}}}
\put(80,-26){\shortstack{\small{$-b$}}}
\put(60,-70){\shortstack{Figure 1.}}
\end{picture}

Let $z_1$ be a coordinate on the first component, $z_2$ be a
coordinate on the second component. Suppose that, the points of
the intersection on the first component have coordinates
$a,-a\in{\mathbb R}$, and on the second $b,-b\in{\mathbb R}$. Let
$$
 P_1=\infty\in\Gamma_1,\
 P_2=\infty\in\Gamma_2,\ r=0\in\Gamma_1,
$$
$$
 Q_1,Q_2,Q_3\in\Gamma_2,\ Q_i\in{\mathbb R},\ \ \gamma\in\Gamma_2,\gamma=i\Gamma\in i{\mathbb R}.
$$
The curve $\Gamma$ has the holomorphic involution
$$
 \sigma :\Gamma\rightarrow\Gamma,\ \sigma(z_1)=-z_1,\ \sigma(z_2)=-z_2.
$$
and antiholomorphic involution
$$
 \mu :\Gamma\rightarrow\Gamma,\ \mu(z_1)=\bar{z}_1,\ \mu(z_2)=\bar{z}_2.
$$
Baker--Akhiezer function $\psi$ on $\Gamma$ is defined by the
functions $\psi_1$ and $\psi_2$ on the components $\Gamma_1$ and
$\Gamma_2$
$$
 \psi_1=e^{ixz_1}f_1(x,y),\ \psi_2=e^{iyz_2} \left(f_2(x,y)+\frac{g_2(x,y)}{z_2-\gamma}\right).
$$
The functions $f_1,f_2$ and $g$ are found from the consistency
conditions
$$
 \psi_1(x,y,a)=\psi_2(x,y,b),\ \psi_1(x,y,-a)=\psi_2(x,y,-b),
$$
and normalization condition
$$
 \psi_1(x,y,0)=d.\
$$
Hence
$$
 f_1=d,\ f_2=\frac{de^{-i(ax+by)}}{2b}(b(e^{2iax}+e^{2iby})+\gamma(-e^{2iax}+e^{2iby})),
$$
$$
 g_2=\frac{de^{-i(ax+by)}}{2b}(e^{2iax}+e^{2iby})(b^2-\gamma^2).
$$
The meromorphic form $\Omega$ is defined by the forms
$$
 \Omega_1=\frac{dz_1}{z_1(z_1^2-a^2)},\
 \Omega_2=\frac{c_1(z_2^2-\gamma^2)dz_2}{(z_2-Q_1)(z_2-Q_2)(z-Q_3)(z_2^2-b^2)},
$$
therefore
$$
 d=\sqrt{\frac{1}{|{\rm Res}_0\Omega_1|}}=a.
$$
Since
$$
 {\rm Res}_a\Omega_1+{\rm Res}_b\Omega_2=0,\  {\rm Res}_{-a}\Omega_1+{\rm
 Res}_{-b}\Omega_2=0
$$
we get
$$
 c_1=-\frac{b(b-Q_1)(b-Q_2)(b-Q_3)}{a^2(b^2-\gamma^2)},
$$
$$
 Q_3=-\frac{b^2(Q_1+Q_2)}{b^2+Q_1Q_2}.
$$
From the condition that the form $\Omega$ has a expansions (6) and
(7) in the neighborhoods $P_1$ and $P_2$ we get $Q_2=Q_1$. Hence
the components of the mapping $\varphi$ have the form
$$
 \varphi_1=\alpha_1F_2(Q_1)=
$$
$$
 \frac{\alpha_1 ae^{-i(ax+(b-Q_1)y)}(e^{2iax}(b+Q_1)(b-i\Gamma)-e^{2iby}(b-Q_1)(b+i\Gamma))}
 {2b(Q_1-i\Gamma)},
$$
$$
 \varphi_2=\alpha_2F_2(Q_2)=
$$
$$
 \frac{\alpha_2 ae^{-i(ax+(b+Q_1)y)}(e^{2iax}(b-Q_1)(b-i\Gamma)-e^{2iby}(b+Q_1)(b+i\Gamma))}{2b(-Q_1-i\Gamma)},
$$
$$
 \varphi_3=\alpha_2F_2(Q_3)=\alpha_3 ae^{-i(ax+\frac{b(b+Q_1)(b+Q_2)}{b^2+Q_1Q_1}y)}\times
$$
$$
 \frac{(-be^{2iax}(b-Q_1)(b-Q_2)(b-i\Gamma)+be^{2iby}(b+Q_1)(b+Q_2)(b+i\Gamma))}
 {2(b^3(Q_1+Q_2+i\Gamma)+ibQ_1Q_2\Gamma)},
$$
where
$$
 \alpha_1=\sqrt{{\rm Res}_{Q_1}\Omega_2}=\sqrt{\frac{b^2(Q_1^2+\Gamma^2)}{2a^2Q_1^2(b^2+\Gamma^2)}},
$$
$$
 \alpha_2=\sqrt{{\rm Res}_{Q_2}\Omega_2}=\sqrt{\frac{b^2(Q_1^2+\Gamma^2)}{2a^2Q_1^2(b^2+\Gamma^2)}},
$$
$$
 \alpha_3=\sqrt{{\rm Res}_{Q_3}\Omega_2}=\sqrt{\frac{\Gamma^2(Q_1^2-b^2)}{a^2Q_1^2(b^2+\Gamma^2)}}.
$$
For $a=b=1,Q_1=2,\gamma=i$ we get
$$
 \varphi_1=\frac{(1+3i)}{8\sqrt{5}}e^{-i(x-y)}(-3ie^{2ix}+e^{2iy}),
$$
$$
 \varphi_2=\frac{e^{-i(x+3y)}}{8\sqrt{5}}((1-3i)e^{2ix}+(9+3i)e^{2iy}),
$$
$$
 \varphi_3=\frac{1}{2}\sqrt{\frac{1}{2}}(\cos(x-y)-\sin(x-y)),
$$
and meanwhile
$$
 e^{2i\beta}=-1.
$$
Induced metric on the image has the form
$$
 ds^2=dx^2+\frac{3}{2}(1+\sin(2(x-y)))dy^2.
$$
The Gaussian curvature of the surface is equal to 1. Hence the
image of the mapping is sphere.

\section*{Acknowledgments}

The author thanks the organizers for the hospitality during the 16th
Osaka City University International Academic Symposium 2008 "Riemann
Surfaces, Harmonic Maps and Visualization" where the work on this
paper was initiated.

\end{document}